\documentstyle[12pt]{article}

\newlength{\abstractwidth}
\renewcommand{\thefootnote}{\fnsymbol{footnote}}
\renewcommand{\thanks}[1]{\footnote{#1}} 
\newcommand{\starttext}{ \setcounter{footnote}{0}
\renewcommand{\thefootnote}{\arabic{footnote}}}

\newcommand{\be}{\begin{equation}}
\newcommand{\bea}{\begin{eqnarray}}
\newcommand{\eea}{\end{eqnarray}} \newcommand{\ee}{\end{equation}}

\def\ba{\begin{eqnarray}}
\def\ea{\end{eqnarray}}

\def\K{{\cal K}}

\def\ra{\rightarrow}

\def\o{\omega}

\def\det{{\rm det}}

\def\z{{\bf z}}

\def\log{\,{\rm log}\,}

\def\o{\omega}

\def\al{\alpha}
\def\b{\beta}

\def\d{\delta}
\def\e{\varepsilon}

\def\l{\lambda}

\def\o{\omega}

\def\t{\theta}
\def\z{\zeta}

\def\O{\Omega}

\def\vp{\varphi}

\def\ti{\tilde}

\def\C{{\bf C}}
\def\P{{\bf P}}

\def\i{\infty}
\def\I{\int}

\def\s{\sum}

\def\ddb{{\partial\bar\partial}}

\def\sub{\subseteq}
\def\ra{\rightarrow}

\def\K{{K\"ahler\ }}

 \def\v{\vskip .1in}

\def\[{{\bf [}}
\def\]{{\bf ]}}

\def\pl{\partial}

\begin{document}
\starttext \baselineskip=15pt \setcounter{footnote}{0}
\newtheorem{theorem}{Theorem}
\newtheorem{lemma}{Lemma}
\newtheorem{definition}{Definition}
\newtheorem{proposition}{Proposition}

\begin{center}
{\Large \bf ON THE SINGULARITIES OF THE PLURICOMPLEX GREEN'S FUNCTION
\footnote{Work supported in part by DMS-07-57372 and DMS-09-05873.
Key words: a priori estimates, K\"ahler manifolds, inequality of Blocki, blow-ups.
AMS Classification: 32W20, 32Q26, 32Q20, 52C44, 53C55.}}
\bigskip\bigskip

{\large  D.H. Phong and Jacob Sturm} \\

\end{center}

\medskip

\begin{abstract}

It is shown that, on a compact K\"ahler manifold with boundary, the singularities of the pluricomplex Green's function with multiple poles can be prescribed to be of the form $\log\sum_{j=1}^n|f_j(z)|^2$ at each pole, where $f_j(z)$ are arbitrary local holomorphic functions with the pole as their only common zero. The proof is a combination of blow-ups and recent a priori estimates for the degenerate complex Monge-Amp\`ere equation, and particularly the $C^1$ estimates away from a divisor.

\end{abstract}

\section{Introduction}
\setcounter{equation}{0}

The Green's function plays a central role in the study of functions of one complex variable or of two real variables. But while its natural generalization to functions of more real variables is as the fundamental solution of the Laplacian, its natural generalization to functions of several complex variables is rather as a fundamental
solution of the complex Monge-Amp\`ere equation. 
For our purposes, we shall consider the following broad definition. 
Let $M$ be an $n$-dimensional compact K\"ahler manifold with smooth boundary $\partial M$, and let $\omega$ be a smooth non-negative closed $(1,1)$-form on $M$.
Let $PSH(M,\omega)$ the space of plurisubharmonic functions with respect to
$\omega$, i.e., $f\in PSH(M,\omega)$ if and only if $f$ is upper semi-continuous, and $\o+{i\over 2}\ddb f\geq 0$ on $X$ in the sense of currents.
Let $\{p_1,\cdots,p_N\}$ be $N$ distinct points on $M$. Then a pluricomplex Green's function $G(z;p_1,\cdots,p_N)$ with poles at the points $p_j$
is a function in $PSH(M,\omega)$, bounded from above on $M$, and
bounded on any compact subset of $\bar M\setminus\{p_1,\cdots,p_N\}$,
which satisfies the equation
\bea
\label{MA}
(\o+{i\over 2}\ddb G)^n=0
\ \ {\rm on}\ M\setminus\{p_1,\cdots,p_N\}
\eea
in the sense of pluripotential theory.
In this paper, we mostly restrict ourselves to the case $\partial M\not=\emptyset$, in which case we also impose the Dirichlet condition
\bea
\label{Dirichlet}
{\rm lim}_{z\to\pl M}G(z;p_1,\cdots,p_N)=0.
\eea

\medskip
We note that the above conditions imply that $(\o+{i\over 2}\ddb G)^n$
is a linear combination of Dirac measures supported at the poles $p_j$. Indeed,
the non-linear expression $(\o+{i\over 2}\ddb G)^n$ is well-defined 
by the Bedford-Taylor construction \cite{BT} as an $(n,n)$ non-negative current away from the unbounded locus of $G$ (i.e., the complement in $M$ of the largest open set where $G$ is locally bounded). Since the unbounded locus consists only of isolated points, and hence is included in compact subsets of a finite union
of Stein neighborhoods, the expression $(\o+{i\over 2}\ddb G)^n$
can again be defined as a non-negative $(n,n)$ current near these points
by the constructions of Demailly \cite{D} and Sibony \cite{Si}.  

\medskip
It has been known for a long-time that the pluricomplex Green's function is not unique, even in the case of a simple pole $p$ \cite{BT}, and that its singularities near $p$ are not unique either. This is in marked contrast with the real Monge-Amp\`ere equation, where the convex solution 
in the sense of Alexandrov on a convex domain $\Omega$ of the equation ${\rm det}\,D^2u=c\,\delta_0, \ \ u_{\pl\Omega}=0$ with $0\in\Omega\subset{\bf R}^n$, $c>0$, is unique. The graph of $u$ is in this case
just an inverted cone, with boundary given by the boundary of $\Omega$, and vertex a point on the $u$ axis, determined by the constant $c$. In particular, the singularities of $u$ at $0$ are determined by $\pl \Omega$ and $c$. This difference between the real and the complex case can be partly attributed to the fact that the condition of convexity in the real case is much more stringent than the condition of plurisubharmonicity.

\medskip
The pluricomplex Green's function has been extensively studied over the years, using many different methods. It is not possible for us to provide a full list of references, but we shall try and indicate along some of the works closest in spirit to the present paper. One notable such example is the work of B. Guan \cite{Gb}, using PDE methods. There, building on the estimates of Yau \cite{Y78} and Caffarelli, Kohn, Nirenberg, and Spruck \cite{CKNS}, he established the existence and $C^{1,1}$ regularity of the pluricomplex Green's function for strongly pseudoconvex domains in ${\bf C}^n$, with prescribed singularity $\log\,|z-p|^2$ at the pole. The $C^{1,1}$ regularity was also obtained by Blocki \cite{B00}. The methods of Guan extend to a singularity of the form $\log\,\sum_{j=1}^n|f_j(z)|^2$, as long as the
$f_j(z)$ are holomorphic functions with only $p$ as their common zero, and are globally defined on the domain.

\medskip
The present paper has two goals. One goal is to establish the existence of pluricomplex Green's functions with singularities at multi poles $p_j$, $1\leq j\leq N$, given by arbitrary local analytic functions. The other is to begin developing a geometric/analytic approach to Monge-Amp\`ere equations with measures on the right-hand side, where the singularities of the solution arise from blow-up constructions. Since blow-ups typically lead to degenerate K\"ahler forms, an essential tool in our approach is the recent existence theorems for the Dirichlet problem for complex Monge-Amp\`ere equations with degenerate background form established in \cite{PS09}.

\bigskip
{\bf Acknowledgements}:
The first author would like to express his profound gratitude to Professor Elias M. Stein for his teachings over the years, beginning with the time when the first author was an undergraduate at Princeton. Both authors would like to thank Professor Stein for the privilege of having worked with him. They would like to dedicate this paper to him, on the occasion of the conference celebrating his 80th birthday.

\section{The Main Results}
\setcounter{equation}{0}

Henceforth, $(M,\omega)$ will be a compact complex manifold of dimension $n$, with non-empty smooth boundary $\pl M$. Let $N$ be an arbitrary positive integer, and let $p_1,\cdots,p_N$ be $N$ distinct (interior) points in $M$. Our main result is the following:

\begin{theorem}
\label{1}
Let $\o$ be a K\"ahler form on $\bar M$, and 
for each $1\leq m\leq N$, let $f_{jm}(z)$, $1\leq j\leq n$, be $n$ holomorphic functions defined in a neighborhood of $p_m$, with $p_m$ as their only common zero in this neighborhood. Then there exists a constant $\e_0>0$ so that for all constants $0<\e_m<\e_0$, $1\leq m\leq N$,
there exists a unique function $G(z;p_1,\cdots, p_N)\in PSH(M,\omega)\cap C^\alpha(\bar M\setminus\{p_1,\cdots,p_N\})$, which satisfies the equation
(\ref{MA}), the boundary condition (\ref{Dirichlet}), and the following asymptotics near each pole $p_m$,
\bea
\label{asymptotics}
G(z;p_1,\cdots,p_N)
=
\e_m\log(\sum_{j=1}^n|f_{jm}(z)|^2)+O(1).
\eea
Here $\alpha$ is any constant satisfying $0<\alpha<1$.
\end{theorem}

\medskip

In general, it is not possible to choose $\e_0$ to be arbitrary.
 In fact, if $\d>0$ is given, then it is easy to construct $(M,\o)$ and $p\in M$, 
and local holomorphic functions $f_1,...,f_n$, for which the maximal $\e_0$ is
less than $\d$. For example, let $M=X\times D$ where $(X,\o)$ is a compact
\K manifold with unit volume, and $D\sub\C$ is the unit disk. Let $p=(x,0)\in X\times \{0\}$ and
choose local coordinates $z_1,...,z_n$ on $X$ centered at $x$. Suppose
$G$ is a Greens function on $M$ with singularity $\e\log (|w|^2+|z_1|^{2k}+\cdots +|z_n|^{2k})$, where $w$ is a coordinate on $\C$, centered at $0\in\C$. Then
$G(0,z_1,...,z_n)\in PSH(X,\o)$ has Lelong number $\e k$ so, by a result
of Demailly \cite{D},
its Monge-Amp\`ere mass is at least $\e k$. Thus $\e \leq {1\over k}$.

\medskip
However, the restriction on $\e_0$ can be removed for
strongly pseudoconvex manifolds, i.e., manifolds $M$ admitting a $C^2$ function $\rho$ with $\pl M=\{\rho=0\}$, and $i\ddb\rho>0$. 
In this case, we have the following solution of the Dirichlet problem:

\medskip

\begin{theorem}
\label{stronglypseudoconvex}
Let  $\omega$ be a non-negative smooth $(1,1)$-form
on $\bar M$. Let $f_{jm}$, $1\leq j\leq n$, $1\leq m\leq N$
be as in the previous theorem.
Assume that $M$ is strongly pseudoconvex.
Then for any function $\vp_b\in C^2(\pl M)$,
and any constant $\e_m>0$, $1\leq m\leq N$,
there exists 
a unique function $G(z;p_1,\cdots, p_N)\in PSH(M,\omega)\cap C^\alpha(\bar M\setminus\{p_1,\cdots,p_N\})$,
which satisfies the equation
(\ref{MA}), the asymptotics (\ref{asymptotics}) near each pole $p_m$,
and the Dirichlet boundary condition
\bea
{\rm lim}_{z\to \pl M}G(z;p_1,\cdots,p_N)=\vp_b.
\eea
Again $\alpha$ is any constant satisfying $0<\alpha<1$.
\end{theorem}

\medskip
Theorem \ref{stronglypseudoconvex} is closest to earlier results of B. Guan \cite{Gb}, Blocki \cite{B00}, and Lempert \cite{L}. The regularity properties of the pluricomplex Green's function are much more precise in these earlier works \cite{Gb}, \cite{B00}, and \cite{L}, but we gain here in generality, including the property that local holomorphic singularities can be assigned arbitrarily. We shall make use of this latter fact
in the construction of new geodesic rays in the space of K\"ahler potentials
(see section \S 5).

\medskip

In general, the case of K\"ahler manifolds with boundary is quite different from the case of K\"ahler manifolds without boundary. Nevertheless,
the idea used in the proof of Theorem 1 may also be used to yield results
for singular Monge-Amp\`ere equations on compact manifolds without boundary.
The simplest example is the following:

\medskip

\begin{theorem}\label{2} Let $(X,\o)$ be a compact \K manifold with unit volume. Let $p\in X$, let
$f\in C^\i(X)$ be positive, and assume $\I_X f\o^n=1$. 
Let $\d_p$ be the Dirac measure concentrated at~$p$.
Then for $\e>0$ sufficiently small, there exists a unique
$\varphi\in PSH(X,\o)\cap C^{1,\al}(X\backslash\{p\})$  satisfying $\varphi = \e\log |z|^2 + C^{1,\al}(X)$ near $z=p$ and
\bea
(\o+{i\over 2}\ddb\varphi)^n\ = \ (1-\e)f\o^n+\e\d_p
\eea
\end{theorem}

\medskip
As Coman-Guedj  \cite{CG} have shown, there are examples of \K manifolds for
which $\e$ must be strictly smaller than one.

\medskip
We observe that the equation $(\o+{i\over 2}\ddb \vp)^n=\mu$
on a compact K\"ahler manifold without boundary has been solved by Berman,
Boucksom, Guedj, and Zeriahi \cite{BEGZ}, when $\mu$
is a measure which does not charge pluripolar sets. The important case of measures $\mu$ which charge pluripolar sets remains open. A general result is that of Ahag et al. \cite{ACCH}, who show that, on a hyperconvex domain in ${\bf C}^n$, a non-negative measure is a complex Monge-Amp\`ere measure if it is dominated by a Monge-Amp\`ere measure. The regularity
and precise singularities of the solutions are however still obscure at this moment. 
Theorem \ref{2} provides another example of the solvability of a Monge-Amp\`ere equation with measures charging pluripolar sets.
In fact, it can be seen from its proof in section \S 4 that more general formulations are possible, with the singularity $p$ replaced by the complex variety
$Z=\{s_1(z)=\cdots=s_k(z)=0\}$, $1\leq k<n$,
where the $s_\alpha(z)$ are sections of a holomorphic vector bundles satisfying some non-degeneracy conditions.

A related and important problem is to provide flexible characterizations of
the unbounded functions for which the Monge-Amp\`ere measure is well-defined
(see e.g. Cegrell \cite{C} and Blocki \cite{B04}).

\section{Proof of Theorems \ref{1} and \ref{stronglypseudoconvex}}
\setcounter{equation}{0}

In this section we give the proof of Theorems \ref{1} and \ref{stronglypseudoconvex}. It will be seen that the argument does not depend
essentially on $N$, so we set $N=1$ and drop the index $m$
to lighten the notation.
The bulk of the work is the proof of Theorem \ref{1},
and we can just indicate at the end the easy modifications for
for Theorem \ref{stronglypseudoconvex}. The key lemmas are the following:

\medskip

\begin{lemma}
\label{lemma1}
Denote by $(M,p,f)$ the data consisting of the K\"ahler manifold $(M,\o)$, the point
$p$, and the given local holomorphic functions $f_j(z)$,
$1\leq\alpha\leq n$ with $p$ as their only common zero. Then there exists a complex manifold $X'=X'(M,p,f)$ and a holomorphic map
$\pi': X'\ra M$ sending $\pl X'$ to $\pl M$, with the following properties:

{\rm (i)} There is a closed, non-negative $(1,1)$-form $\Omega'$ on $X'$
an effective divisor $E'$, and an $\e>0$ with
\bea
\Omega'-\e{ i\over 2}\ddb \log h_{E'}>0
\eea
for some smooth metric $h_{E'}$ on $O(-E')$.

{\rm (ii)} The restriction $\pi=\pi'|_{\bar X'\backslash E'}$ is a biholomorphism $\pi:\bar X'\setminus E'
\to \bar M\setminus p$ with
\bea
\label{Omegaprime}
\pi_*\Omega'
=
\o+\e {i\over 2}\ddb (\psi(z)\,\log \sum_{j=1}^n|f_j(z)|^2 + 1-\psi(z))
\eea
where $\psi(z)$ is a function which is $1$ in a neighborhood of $p$,
and which is compactly supported in another such neighborhood.
\end{lemma}

\medskip

\begin{lemma}
\label{lemma2}
Let $X'$ be a complex manifold with smooth boundary
of dimension $n$, equipped with a non-negative closed form $\Omega'$, with $\Omega'-\e{i\over 2}\ddb\log h_{E'}>0$ for some effective divisor $E'$ supported away from $\pl M$
some smooth metric $h_{E'}$ on $O(-E')$, and some $\e>0$. Then there exists a unique function $\Phi\in PSH(X',\Omega')\cap L^\infty(X')\cap C^\alpha(\bar X'\setminus E')$ which solves the Dirichlet problem
\bea
\label{HCMA}
(\Omega'+{i\over 2}\ddb \Phi)^n=0
\ \ {\rm on}\ X',
\qquad
\Phi_{\pl X'}=0.
\eea
Here $\alpha$ is any constant satisfying $0<\alpha<1$.
\end{lemma}

\medskip
The theorem follows readily from the two lemmas. Let $\Phi$ be the function
given by Lemma \ref{lemma2} applied to the complex manifold $X'=X'(M,p,f)$ and the non-negative form
$\Omega'$ of (\ref{Omegaprime}). 
Then  $(\Omega'+{i\over 2}\ddb\Phi)^n=0$ on $X'$.
Set $\vp=\Phi\circ \pi$. 
Since $\pi$ is a biholomorphism between $X'\setminus E'$ and $M\setminus p$,
this implies $(\pi_*\O'+{i\over 2}\ddb\phi)^n$ on $M\backslash p$, i.e.
\bea
0=
(\omega+{i\over 2}\ddb (\e[\psi(z)\,\log \sum_{j=1}^n|f_j(z)|^2+
(1-\psi(z))]+\vp(z)))^n.
\eea
We can now set
\bea
G(z;p)=\e [\psi(z)\,\log \sum_{j=1}^n|f_j(z)|^2+(1-\psi(z))]+\vp(z)-\e.
\eea
Clearly it satisfies the equation $(\o+{i\over 2}\ddb G)^n=0$
on $M\setminus p$, vanishes on $\pl M$, and has the desired asymptotics near $p$. It is $\omega$-plurisubharmonic on $M\setminus p$, and bounded from above. Thus it extends to an $\omega$-plurisubharmonic function on $M$. This shows that $G(z;p)$ satisfies all the desired properties, and the existence part of the theorem is proved.
\v
To prove the uniqueness, let $\ti G\in PSH(M,\o)$ satisfy (\ref{MA}),
(\ref{Dirichlet}) and (\ref{asymptotics}). Let
\bea
\ti\Phi = \ti G-\e [\psi(z)\,\log \sum_{j=1}^n|f_j(z)|^2+(1-\psi(z))]
+\e
\eea
Then $\ti\Phi\in PSH(\O', X'\backslash E')\cap C^\al(X'\backslash E')$. Since it is bounded, it extends to a function, which by abuse of notation,
will be denoted
$\ti\Phi\in PSH(\O', X')\cap C^\al(X'\backslash E')\cap L^\i(X)$.
Moreover, $(\O'+{i\over 2}\ddb\Phi')^n=0$. This is certainly true on
$X'\backslash E'$ and, since $\Phi'$ is bounded, is true on all of $X'$.
By Lemma 2, we have $\ti\Phi=\Phi$ and hence, $\ti G=G$.

\medskip

Thus it suffices to establish Lemma \ref{lemma1} and Lemma~\ref{lemma2}.
We begin with the proof of Lemma~\ref{lemma1}, which we break into several steps. Some of these are well-known, but we did not find the version that we needed in the literature, so we have provided a complete derivation. 

\subsection{Blow-ups of complex manifolds}

We start by recalling the construction 
of the blow-up of a complex manifold along a smooth
submanifold and some of its basic
properties. 

\medskip

Let $W$ be a complex manifold of dimension $n$, and let $Z\subset W$ be a submanifold of dimension $d<n$. Let $N=TW_{\vert_Z}/TZ$ be the normal bundle of $Z$ and let $E={\bf P}(N)$. Let
\bea
W'=(W\setminus Z)\cup E
\eea
and define $\pi: W'\ra W$ by extending the map $\pi:E\ra Z$ to be the
identity map on $W\backslash Z$. The set $W'$ has a natural complex
structure for which $\pi: W'\ra W$ is a holomorphic map:
The complex structure on $W'$ is defined as follows: 

\medskip
First, we require that
$W\backslash Z\sub W'$ is an open set and the inclusion $W\backslash Z\ra W'$
is holomorphic.
\v
Next, we let $U_\al$ be a collection of coordinate balls in $W$ such that $Z\sub \cup_\al U_\al$ and such that 
$U_\al\cap Z = \{(x^\al,y^\al)\in U_\al: y^\al=0\}$. On $U_\al\cap U_\b$
we have $z^\b=\phi^\b_\al(z^\al)$ with $\phi^\b_\al: U_{\al\b}\ra U_{\b\al}$ a biholomorphic
function between open subsets of $\C^n$. The derivative $D\phi^\b_\al
=\left(\matrix{D_{11} & D_{12}\cr D_{21} & D_{22}\cr}\right)$
is an invertible $n\times n$ matrix, where $D_{11}$ has size $d\times d$,
and $D_{22}$ has size $(n-d)\times (n-d)$. 
\v
Note that $D_{22}(p)$ is invertible
if $p\in U_\al\cap Z=\{z^\al\in U_\al: y^\al=0\}$. In fact, $D_{22}$ is the isomorphism of the normal bundle of
$\{y^\al=0\}\sub U_{\al\b}$ to 
the normal bundle of
$\{y^\b=0\}\sub U_{\b\al}$ induced by the biholomorphic map $\phi^\b_\al$.
Let us spell out this point in more detail. We can write
\bea
y^\b_i= A^j_iy^\al_j, \ \ 1\leq i,j\leq n-d,
\eea
for certain (non-unique) holomorphic functions $A^l_k$ on $U_{\al\b}$.
This follows from the fact that $\phi^\b_\al$ takes the $x^\al$-axis (i.e.
the set $\{y^\al=0\}$) to
the $x^\b$-axis
(i.e.
the set $\{y^\b=0\}$). Then
$D_{22}(p) = (A^j_i(p))$ is an invertible matrix.

\v
Define a complex manifold $U_\al'$ as follows: 
\bea
U_\al'=\{(z,t): z\in U_\al, t\in\P^{n-d-1}: y_it_j=t_iy_j, 
1\leq i,j\leq n-d\}
\eea
Let $f_\al: U_\al'\ra \pi^{-1}(U_\al)\sub W'$ be the bijective map
\bea
(z,t)\mapsto \cases{ z & if $y\not=0$ \cr
t_1{\pl\over \pl y_1}+\cdots + t_{n-d}{\pl\over y_{n-d}} & if $y=0$
}
\eea
We wish to use the maps $f_\al$ to give $W'$ a complex structure.
To do this, we must show that the change of coordinate maps
$f^\b_\al=f_\b^{-1}\circ f_\al: U_{\al\b}'\ra U_{\b\al}'$ are holomorphic.
\v
We proceed as follows: $f^\b_\al(z^\al,t^\al)=(z^\b,t^\b)$
where $z^\b=\phi^\b_\al(z^\al)$ and
\bea 
{t^\b_j\over t^\b_i}={y^\b_j\over y^\b_i}= {A^k_jy^\al_k\over A^k_iy^\al_k}=
{A^k_jt^\al_k\over A^k_it^\al_k}
\eea
If $y^\b_i\not= 0$ for some $i$, then the first equality implies that
$t^\b_i\not=0$ so we can take $t^\b_i=1$ and $t^\b_j={y^\b_j\over y^\b_i}$,
which is holomorphic.
\v
If $y^\b_i=0$ for all $i$ then we make use of the fact that $D_{22}$ is
invertible on $Z$ so there exists $i$ such that $A^k_it^\al_k\not=0$.
We take $t^\b_i=1$ and $t^\b_j={A^k_jy^\al_k\over A^k_iy^\al_k}$
which is holomorphic. Thus we see that $W'$ is a complex manifold.

\v
Now let $\pi:W'\ra W$ be as above. Let $p\in E$ and $q=\pi(p)$. Then the discussion
above shows that there exists a coordinate neighborhood $\O$ of $p\in W'$ and
coordinates $(\z_0,....,\z_d, \t_1,...,\t_{n-d-1})$ centered at
$p$ with the following properties:
\v
1) $\z_i=z_i\circ\pi$ for some coordinate
functions $z_j$ on $W$
\v
2) $E\cap\O=\{\z_0=0\}$.
\v
3) $(\t_1|_E,...,\t_{n-d-1}|_E)$ is a set of local coordinates
of $\pi^{-1}(p)$ centered at $p\in E$. 
\v
4) If $p\in U_\al'$ then $\z_0| y^\al_j$ for all $j$ and
$y_{j_0}^\al/\z_0$ is nowhere vanishing for some $j_0$.
\v
The last condition says that $E\cap\O=\{y^\al_{j_0}=0\}$.

\v

This blow-up process can be iterated: If $Z'\sub W'$ is a smooth subvariety then
we can construct $W''=BL(Z',W')$ and we have maps $W''\ra W'\ra W$. We
say that $W''$ is an iterated blow-up. 
\v
If $\pi:W'\ra W$ is an iterated blow-up, the exceptional divisor is by definition
the smallest effective divisor $E\sub W'$ such that $W'\backslash E \ra W\backslash \pi(E)$
is an isomorphism.

\subsection{Analytic spaces}

\medskip
Let $X$ be a set. We say that $X$ is an analytic space if there is a complex manifold $W$ (called an ambient manifold)
with $X\subset W$ and satisfying the following property: for every $p\in X$,
there is an open set $p\in U\subset W$ and functions $f_1,\cdots,f_r:U\to {\bf C}$
such that 
\bea
U\cap X=\{w\in W;\, f_1(w)=\cdots=f_r(w)=0\}.
\eea
We denote by $X_{reg}$ the subset of $X$ where $X$ is locally a complex manifold, and by $X_{sing}$ its complement in $X$. A K\"ahler metric on $X$ is by definition a K\"ahler metric on $X_{reg}$ which is the restriction of a K\"ahler metric on an ambient manifold $W$.

\medskip
\begin{lemma}
\label{lemmaX}
Let $(M,p,f)$ be a data consisting of a complex manifold $M$,
a point $p\in M$, and local holomorphic functions $f_1(z),\cdots,f_n(z)$
defined in a neighborhood of $p$, and with $p$ as their only common zero. 
Then we can associate to this data an analytic space
$X=X(M,p,f)$ with the following properties:

{\rm (i)} There is a biholomorphism between $X\setminus X_0$ and $M\setminus p$, for some subset $X_0$ of $X$ which is biholomorphically equivalent to ${\bf CP}^{n-1}$.

{\rm (ii)} Let $\psi(z)$ be a cut-off function which is $1$ in a coordinate chart around
$p$ in $M$, and $0$ outside another such chart. Then for all $\delta$ small enough, the pull-back to
$X\setminus X_0$ of the form $\o_\delta$ defined as
\bea
\label{omegadelta}
\o_\delta=\o+\delta{i\over 2}\ddb (\psi(z)\log \sum_{j=1}^n|f_j(z)|^2+1-\psi(z))
\eea
defined on $M\setminus p$ extends to a K\"ahler form on $X$.

\end{lemma}

\medskip
{\it Proof of Lemma \ref{lemmaX}}: Given the data $(M,p,f)$, let $U$ be a coordinate neighborhood centered at $p$ in $M$,
define a space $V$ by
\bea
V=\{(z_1,...,z_n),(y_1,...,y_n)\in U\times {\bf CP}^{n-1}: y_if_j(z)=y_jf_i(z)\}
\eea
and let $\pi: V\ra U$ be the projection. We then define 
$X=X(M,p,f)$ by
\bea
X= (M\backslash \{p\}\cup V)/\sim
\eea 
where, for $m\in M$ and $v\in V$
we say $m\sim v$ if $\pi(v)=m$. Note that the fiber $X_0$
of $V$ above the point $p$ is the entire projective space ${\bf CP}^{n-1}$.
We claim that $X$ is an analytic space.
\v
To show that $X$ is an analytic space, we must 
find a complex manifold $W$ such that $X\sub W$ and such that $X$ is locally defined
by the simultaneous vanishing of a finite collection of holomorphic functions.
\v
Let $B\sub\C^n$ be a small open ball centered at the origin and
let $Z\sub U\times \C^n$ be the smooth manifold
\bea
Z=\{(z,\xi): z\in U, \xi\in B: f_1(z)-\xi_1=\cdots=f_n(z)-\xi_n=0\}
\eea
If $B$ is sufficiently small, then the image of the map $Z\ra U$
is compactly supported in $U$. Thus $Z\sub M\times B$
is a smooth submanifold whose image, when projected to $M$, lies in
a relatively compact subset of $U$. Finally, let $W=BL(Z, M\times B)$.
Thus $W$ is locally defined by
\bea
W=\{(z,\xi,y)\in U\times B\times {\bf CP}^{n-1}: y_i(f_j(z)-\xi_j)=y_j(f_i(z)-\xi_i)\}
\eea
Then $W$ is a smooth manifold and $X\sub W$ is  defined by $\xi_1=\cdots=\xi_n=0$. This shows $X$ is an analytic space.
\v
We can define a \K metric on $X$ as follows. Let $\o$ be a \K metric on $M$.
Extend $\o$ to a \K metric on $M\times B$. Choose $\psi
\in C^\i(U)$ with 
the property that $\psi$ equals one in a neighborhood $p$ and 
${\rm support}(\psi)\sub U$
is compact. Then the composition of the map
$U\times B\ra U$ with $\psi$ is a smooth function on $U\times B$ which,
by abuse of notation, is again denoted $\psi$. By ${\rm support}(\psi)\sub U$
to be a sufficiently large compact set, we may assume $\psi=1$ on $Z$.
Let
\be\label{metric} \o_\d = \o + i{\d\over 2}\ddb[(1-\psi)+\psi\log  ( \s |f_j(z)-\xi_j|^2)\
\ee
on $W\backslash E = (M\times B)\backslash Z$ and
let
\be\label{metric1} \o_\d\ = \ \o\ + \ \d({i\over 2}\ddb\log |y|^2)\ = \ \o + \d\o_{\rm FS}
\ee
on the open neighborhood $\{\psi=1\}^o$ of $E$   (i.e., the interior
of the closed set  $\{\psi=1 \}$). Here $\o_{\rm FS}$
is the Fubini-Study metric on $\P{\bf C}^{n-1}$.
The two definitions are consistent, and define
$\o_\d$, a smooth $(1,1)$ form on $W$. Since $\o+\d\o_{\rm FS}>0$ in a
fixed (independent of $\d$) neighborhood of $E$, we find that $\o_\d>0$
on all of $W$
for $\d$ sufficiently small. Thus $\o_\d$
is a \K metric on $W$ and its restriction to $X$ is a \K metric on $X$.
The proof of Lemma~\ref{lemmaX} is complete.

\medskip
For general functions $f_1(z),\cdots,f_n(z)$, the space $X(M,p,f)$ is not smooth. In order to obtain a smooth manifold, we use Hironaka's theorem on resolution of singularities. One version of this theorem is the following:

\medskip
\begin{theorem}
\label{lemmaHironaka}
Let $W$ be a complex manifold and let $X\subset W$ be a complex analytic space. Then there exists an iterated blow up $W'$ of $W$, with corresponding holomorphic map
$\pi:W'\to W$ with the following property: let $E$ be the exceptional divisor, and set
\bea
X'=\overline{\pi^{-1}(X)\setminus E}.
\eea
Then $X'$ is a smooth manifold and the map $\pi:X'\to X$ is surjective. Moreover,
\bea
E'=E\cap X'=\pi^{-1}(X_{sing})
\eea
is a divisor with normal crossings, and $\pi:X'\setminus E'\to X_{reg}$ is an isomorphism.
\end{theorem}

\subsection{Metrics on blow-ups}

The following is the key property of blow-ups that we need:

\begin{lemma}
\label{lemmaKahler}
Let $W$ be a smooth complex manifold and $Z\subset W$ a smooth submanifold. Let $\pi:W'\to W$ be the blow-up of $W$ with center $Z$, and let $E\subset W$ be the exceptional divisor. If $\o$ is any K\"ahler metric on $W$, then there exists a hermitian metric $h_E$ on $O(-E)$ so that, for any compact subset $K$ of $W'$,
there exists $\e_K>0$ with
\bea
\pi^*\o-\e{i\over 2}\ddb \log h_E>0
\eea
for all $0<\e<\e_K$.
\end{lemma}

\medskip
\v
{\it Proof of Lemma \ref{lemmaKahler}:}
Let $\{U_\al\}\sub W'$ be a locally finite collection of coordinate
neighborhoods which cover $Z$, and let $z^\al=(x^\al,y^\al)$
be coordinates on $U_\al$ with the property:
$$ U_\al\cap Z\ = \ \{y^\al=0\}
$$
  Choose $\psi_\al\in C^\i_c(U_\al)$
such that $0\leq \psi_\al\leq 1$ and $\psi=\s_{\al=1}^r \psi_\al =1$ on $Z$. 

\v

Let $f$ be a section of 
$O(-E)$ over an open set $\O$. This means that $f$ is a holomorphic
function on $\O$ which vanishes on $E\cap \O$. 
Let
\bea
|f|^2_{h_E} = {|f|^2\
\over (1-\psi)+\s_\al\psi_\al( |y^\al_1|^2+\cdots + |y^\al_{n-d}|^2)}
\eea
Observe that
this makes sense: We may assume that $E\cap\O\not=\emptyset$ and that
$\O$ is a small open set such that $\psi\circ\pi=1$ on $\O$.
Choose coordinates $(\z,\t)$ as in section 13.1. Thus
$\z=(\z_0,...,\z_d)$ and $\t=(\t_1,...,\t_{n-d-1})$. Since $f$ 
vanishes on $E$ we have $f=g\z_0$ for some holomorphic function $g$ and
\bea
|f|^2_{h_E} = {|f|^2\
\over \s_\al\psi_\al( |y^\al_1|^2+\cdots + |y^\al_{n-d}|^2)} = 
{|g|^2\over 
\s_\al\psi_\al( |y^\al_1/\z_0|^2+\cdots + |y^\al_{n-d}/\z_0|^2)
}
\eea
This shows that $h_E$ is a well defined smooth metric on $O(-E)$.
\v
Next we claim that $(\pi^*\o-{i\over 2}\e\ddb\log h_E)(p)>0$ for all $p\in E$ and
 sufficiently small $\e>0$. To see this, fix $p\in E$. 
 If $F(\z,\t)$ is smooth in a neighborhood of $\{p\}$ then
\bea
\pl_i\pl_{\bar j}F=
\left(
\matrix{ A & B\cr
{}^t\bar B & D\cr
}
\right)
\eea
where $A$ has size $(d+1)\times(d+1)$, $B$ has size $(d+1)\times (n-d-1)$ and
$D$ has size $(n-d-1)\times (n-d-1)$. We see that in these coordinates,
\bea
\o(p)  =
\left(
\matrix{ A & 0\cr
0 & 0\cr
}
\right)
\eea
with $A>0$. Now we write
\bea
-i\ddb\log h_{E}(p)\ = \ 
\left(
\matrix{ X & Y\cr
{}^t\bar Y & D\cr
}
\right)
\eea
Since $\psi$ is independent of $\t_j$ we have
\bea
D_{i\bar j}\ = \ {\pl^2\over \pl\t_j\pl \bar \t_k}\log\s_\al\psi_\al(p)( |y^\al_1/\z_0|^2+\cdots + |y^\al_{n-d}/\z_0|^2)
\eea
Finally we observe that $\s_{i,j=1}^{n-d-1}D_{i\bar j}d\t_i\wedge d\t_{\bar j}$ is the pullback of a  Fubini-Study metric
with respect to a holomorphic map whose deriviative has maximal rank.
Thus $D>0$.
\v
The claim now follows from the following  linear algebra fact:

\begin{lemma} 
\label{linearalgebra}
Let $A,X$ be $(d+1)\times (d+1)$ hermitian matrices
and  $D$  an $(n-d-1)\times (n-d-1)$ hermitian matrix. Let
$Y$ be a $(d+1)\times (n-d-1)$ matrix. Assume $A>0$ and $D>0$. Then
for $\l>0$ sufficiently large we have
\bea
M(\l)\ = \ \l \left(
\matrix{ A & 0\cr
0 & 0\cr
}
\right)
+\left(
\matrix{ X & Y\cr
{}^t\bar Y & D\cr
}
\right)\ > \ 0
\eea
\end{lemma}
{\it Proof.} We may assume $A$ to be diagonal with positive diagonal entries
$a_0,...,a_{d}$. Then $\det(M(\l)$ is a polynomial of degree $d$, with real coefficients in $\l$ whose leading
coefficient is $a_0\cdots a_d>0$. Thus $\det\,M(\l)>0$ for $\l>0$ sufficiently
large.
\v
The same argument shows that the determinants of all the square submatrices
of $M(\l)$ which are situated in the lower right corner of $M(\l)$, also have
positive determinants for $\l$ sufficiently large. This proves Lemma
\ref{linearalgebra}.
\v
Since $(\pi^*\o-{i\over 2}\e\ddb\log h_E)>0$ on $Z$, we see that there is an open neighborhood
$U$ of $Z\cap K$ on which $(\pi^*\o-{i\over 2}\e\ddb\log h_E)>0$ for all sufficiently small $\e>0$. Choosing
$\e>0$ sufficiently small, we have $(\pi^*\o-{i\over 2}\e\ddb\log h_E)>0$ on all of $K$.
This proves Lemma \ref{lemmaKahler}.

\medskip
The resolution of singularities in Hironaka's theorem will usually require an iteration of blow-ups. Thus we need the following generalization of Lemma \ref{lemmaKahler}:

\begin{lemma}
\label{lemmaKahleriterated}
Let $W$ be a compact complex manifold and let $\pi:W'\to W$
be an iterated blow-up of~$W$. Let $E\subset W'$ be the exceptional divisor. If $\o$ is any K\"ahler metric on $W$, then there exists an effective divisor
$E'$ on $W'$
supported on $E$, and a hermitian metric $h_{E'}$ on $O(-E')$ so that
there exists $\e>0$ with
\bea
\pi^*\o-\e{i\over 2}\ddb \log h_{E'}>0.
\eea
\end{lemma}

\medskip
{\it Proof of Lemma \ref{lemmaKahleriterated}:} Let $\pi:W'\to W$ be the composition of two blow-ups, $\pi=\pi_2\circ \pi_1$, with $\pi_2:W'=W_2\to W_1$,
$\pi_1=W_1\to W$. Apply Lemma \ref{lemmaKahler} to $\pi_1$ and the K\"ahler metric $\o$ on $W$. If $E_1$ is the exceptional divisor of $\pi_1$,
we obtain a metric
$h_{E_1}$ in on the line bundle $O(-E_1)$
on $W_1$ with $(\pi_1)^*\o-\e_1{i\over 2}\ddb \log h_{E_1}>0$ on $W_1$. 
Apply next Lemma \ref{lemmaKahler} to $\pi_2$ and the K\"ahler metric
$(\pi_1)^*\o-\e_1{i\over 2}\ddb \log h_{E_1}$ on $W_1$. We obtain then a metric $h_{E_2}$ on $O(-E_2)$, where $E_2$ is the exceptional divisor of $\pi_2$, with
$\pi_2^*(\pi_1^*\o-\e_1{i\over 2}\ddb \log h_{E_1})-\e_2{i\over 2}\ddb \log h_{E_2}>0$
on $W_2$.  We can take $\e_1={1\over n_1}$ and $\e_2={1\over n_1n2}$ for
$n_1$ and $n_2$ large enough integers. We can then write
\bea
\pi_2^*(\pi_1^*\o-\e_1{i\over 2}\ddb \log h_{E_1})-\e_2{i\over 2}\ddb \log h_{E_2}
=
\pi^*\o-{1\over n_1n_2}
{i\over 2}
\ddb
\log h_{n_2(\pi_2)^*E_1}h_{E_2}
\eea
so the lemma holds in this case with the line bundle given by
$O(-E_2)\otimes \pi_2^*O (-n_2E_1)$. Clearly the argument extends to any finite number
of blow-ups, and the lemma is proved.

\bigskip

The preceding lemma can be extended to the case of complex analytic sets:

\begin{lemma}
\label{lemmaKahleranalytic}
Let $X$ be a complex analytic set, and $\o$ a K\"ahler metric on $X$.
Let $\pi:X'\to X$ be a resolution of singularities and $E\subset X'$ the exceptional divisor. Then there is a divisor $E'$ on $X'$, whose support is contained in $E$, a hermitian metric $h_{E'}$ on the line bundle $O(-E')$, and an $\e>0$ 
such that
$\pi^*\o-\e{i\over 2}\ddb\log h_{E'}>0$.
\end{lemma}

\medskip

{\it Proof of Lemma \ref{lemmaKahleranalytic}}: By definition of a K\"ahler metric on $X$, the metric $\o$ extends to a K\"ahler metric on an ambient space $W$ of $X$.
By definition of resolution of singularities, the map $\pi$ extends to a map $\pi:W'\to W$, with $X'\subset W'$, and $W'$ an iterated blow-up. The metric on the line bundle $O(-E')$ on $W'$ obtained from Lemma \ref{lemmaKahleriterated} applied to $W$ and the K\"ahler form $\o$
restricts to a metric on the line bundle $O(-E')$ on $X'$ with the desired property.

\bigskip
We can now give the proof of Lemma \ref{lemma1}: We apply Lemma \ref{lemmaX}, to obtain the analytic space $X=X(M,p,f)$ and the K\"ahler form 
$\o_\delta$ with the properties stated there. The space $X$ is only an analytic set, but we can apply
Theorem \ref{lemmaHironaka} to obtain a resolution of singularities $\pi:X'\to X$. We can then apply Lemma \ref{lemmaKahleranalytic} to obtain a line bundle 
$O(-E')$ on $X'$
with $\pi^*(\o_\delta)-\e{i\over 2}\ddb \log h_{E'}>0$ on $X'$.  
Set
\bea
\Omega'=\pi^*\o_\delta.
\eea
Since the resolution of singularities is a biholomorphism of $X'\setminus E'$ to $X_{reg}$,
and $X_{reg}$ contains (a biholomorphic image of) $M\setminus p$,
the K\"ahler form $\Omega'=\pi^*(\o_\delta)$
retains the same expression (\ref{omegadelta}) on $M\setminus p$ , and is hence given by
the expression (\ref{Omegaprime}). The proof of Lemma \ref{lemma1} is complete.

\bigskip
It remains only to establish Lemma \ref{lemma2}. In the special case when the background form $\Omega'$ is actually strictly positive,
this has been proved by X.X. Chen \cite{C00} and Blocki
\cite{B09}. But for our purposes, it is essential to allow degeneracies in $\Omega'$,
as such degeneracies arise due to blow-ups.
The full Lemma \ref{lemma2}, allowing for degeneracies,
is actually the main result of
\cite{PS09}, stated there as Theorem 2. The desired solution of the homogeneous complex Monge-Amp\`ere equation (\ref{HCMA}) is obtained as a $C^{\alpha}$ limit
on compact subsets of $X'\setminus E$ of solutions of elliptic equations where the right hand side tends to $0$. The key estimate is the following pointwise
$C^1$ estimate (\cite{PS09}, Theorem 1)
\bea
\label{C1}
|\nabla \vp(z)|
\leq C_1 {\rm exp}(C_2\vp(z))
\eea
for the solutions of the Dirichlet problem for the equation
\bea
\label{HCMAa}
(\o+{i\over 2}\ddb \vp)^n= F(z,\vp)\o^n
\eea
where $\o$ is a K\"ahler form on $X'$.
Here $C_1$ and $C_2$ are strictly positive constants which depend only on upper bounds for ${\rm inf}_{X'}\vp$, ${\rm sup}_{X'\times [{\rm inf}_{X'}\vp,\infty)}F$,
${\rm sup}_{X'\times [{\rm inf}_{X'}\vp,\infty)}(|\nabla F^{1\over n}|+|\pl_\vp F{1\over n}|)$, $\|\vp\|_{C^1(\pl X')}$ and a lower bound for the holomorphic bisectional curvature of $\o$. The point of the estimate (\ref{C1}) is that it does not require an upper bound for ${\rm sup}_{X'}\vp$.
It is used to obtain
the existence of convergent subsequences in $C^{\alpha}(X'\setminus E')$ for any $0<\alpha<1$ of solutions $\vp$
of the equation (\ref{HCMAa}) which may tend to $+\infty$ along a divisor $E'$. 
The proof of (\ref{C1}) makes essential use of a differential inequality for solutions of complex Monge-Amp\`ere equations due to Blocki \cite{B09}. We refer to \cite{PS09} for the complete details. The proof
of the main theorem is complete.

\subsection{The case of strongly pseudoconvex manifolds}

We give now the proof of Theorem \ref{stronglypseudoconvex}.

\medskip
First we observe that the theorem can be reduced to the case of boundary value $0$ and $\o$ strictly positive. Indeed, for any boundary value $\vp$, we can pick an extension $\hat\vp$
of $\vp$ to $\bar M$ with the property that $i\ddb\hat\vp>0$. This can be done by choosing any extension, and adding a large positive multiple of $i\ddb \rho$.
Next, set $\hat\o=\o+{i\over 2}\ddb \hat\vp$, and $\hat G=G-\hat\vp$. The equation
(\ref{MA}) can be rewritten as
\bea
(\o+{i\over 2}\ddb G)^n=
(\hat\o+{i\over 2}\ddb \hat G)^n=0\ \ {\rm on}\ \ M\setminus\{p_1,\cdots,p_N\}.
\eea
So if we can solve this equation for $\hat G\in PSH(M,\o)$ with boundary value $0$,
then the function $G=\hat G+\hat\vp\in PSH(M,\hat\o)$ is a solution of the original problem.

Next, assume that the boundary value is $0$ and the form $\o$ is strictly positive. We note that the restriction to $\e_j$ small in the proof of Theorem \ref{1} is just due to the requirement that the form $\o_\delta$ of (\ref{omegadelta}) be strictly positive. But in the present case, for any $\delta$, it suffices to replace the form
$\o_\delta$ by the form $\o_\delta+A(\delta)i\ddb \rho$ with $A$ large enough,
in order to obtain a form which is strictly positive. The rest of the proof applies verbatim. The proof of Theorem \ref{stronglypseudoconvex} is complete.

\section{Proof of Theorem \ref{2}}
\setcounter{equation}{0}

Let $X'=BL(X,p)$, the blow up of the point $p$, and let 
$\pi: X'\ra X$ be the projection map. Choosing, as before,
a cut-off function $\psi$ which is supported in a neighborhood
of the point $p$. Then for $\e$ sufficiently small,
\bea
\o_\e = \o + {i\over 2}\e\ddb(\psi\log |z|^2+(1-\psi))
\eea
extends to a smooth \K metric $\o'$ on the smooth manifold $X'$. 
Consider the equation on $X'$
\bea
(\o'+{i\over 2}\ddb\varphi')^n = cf(\pi^*\o)^n
\eea
for a function $\vp'\in PSH(X',\o')$ with $c$ a normalization constant
so that both sides have the same total volume. This equation can be rewritten as
\bea
\label{MA1}
(\o'+{i\over 2}\ddb\varphi')^n = cF\,(\o')^n
\eea
where $F\equiv f{(\pi^*\o)^n\over(\o')^n}$ is a smooth non-negative function.
A careful examination of Yau's treatment \cite{Y78}
of equations of the form (\ref{MA1})
shows that a priori upper bounds for $\|\vp'\|_{C^0(X')}$ and for
$\|\Delta_{\o'}\vp'\|_{C^0(X')}$ do not require a strictly positive lower bound for $F$.
Thus the equation (\ref{MA1}) admits a generalized solution
$\vp\in PSH(X',\o')\cap C^{1,\alpha}(X')$ for any $0<\alpha<1$.
Restricting to $X'\backslash E$ we get
\bea
(\o + {i\over 2}\ddb(\e(\psi\log |z|^2+1-\psi)+\varphi'))^n = cf\,\o^n
\ \ \hbox{on $X\backslash \{p\}$}
\eea
Thus, if we let $\varphi = \e\psi \log |z|^2+(1-\psi)+\varphi'$ we get
$$ (\o + {i\over 2}\ddb\varphi)^n\ = \ \e\d_p + cf\o^n
$$
which implies $c=1-\e$. The proof of Theorem \ref{2} is complete.
\v
Clearly the proof extends to the cases of local singularities $f_j(z)$,
when the analytic set $X=X(M,p,f)$ is a smooth manifold. It is easy to
formulate conditions on the $f_j(z)$ which would insure this property,
but we leave this to the interested reader.

\section{Geodesics in the space of K\"ahler potentials}
\setcounter{equation}{0}

Let $(X,\o_0)$ be a compact K\"ahler manifold without boundary. A well-known conjecture of Yau \cite{Y93} is that the existence of a K\"ahler form in the class
$[\o_0]$  with constant scalar curvature
should be equivalent to the stability of $(X,[\o_0])$
in the sense of geometric invariant theory. Suitable notions of stability have been
proposed by Tian \cite{T97} and Donaldson \cite{D98, D02} (see also \cite{PS06}
\cite{PSSW} for some other notions of stability, and \cite{PS09a} for a survey).
In particular, in \cite{D98}, Donaldson introduces a notion of stability based on the behavior of the $K$-energy functional of Mabuchi near
infinity along geodesic rays in the space of K\"ahler potentials. Such rays have been constructed from test configurations (see \cite{PS06}\cite{PS07}\cite{CT}\cite{SZ08}\cite{RWN},
and also \cite{AT} in the analytic category, using the Cauchy-Kowalevska theorem).
Here we illustrate Theorem \ref{1} by constructing certain new rays, exploiting the fact that local singularities can be prescribed near infinity.

\medskip
More precisely, the space ${\cal K}$ of K\"ahler potentials is defined by
\bea
{\cal K}=\{\vp\in C^\infty(X);\o_\vp\equiv \o_0+{i\over 2}\ddb\vp>0\}.
\eea
It carries a natural Riemannian structure defined by the $L^2$ norm
on $T_\vp({\cal K})$ with respect to the volume form $\o_\vp^n$.
A path $(-T,0]\ni t\to \vp(\cdot,t)$ is a geodesic if and only if it satisfies the equation
\bea
\label{geodesic}
\ddot\vp-g_{\vp}^{j\bar k}\pl_j\dot\vp\pl_{\bar k}\dot\vp=0.
\eea
where $g_{\vp j\bar k}$ is the metric corresponding to the K\"ahler form $\o_{\vp}$. A key observation due to Donaldson \cite{D98} and Semmes \cite{Se} is that this equation is equivalent to the homogeneous complex Monge-Amp\`ere equation
\bea
\label{HCMA2}
(\pi^*\o_0+{i\over 2}\ddb\Phi)^{n+1}=0
\eea
on the manifold $M=X\times \{e^{-T}<|w|<1\}$, for the function $\Phi$ defined by
\bea
\Phi(z,w)=\vp(z, \log |w|),
\eea
and where $\pl$ is now with respect to both $z$ and $w$. The end points of the geodesic paths in ${\cal K}$ correspond to Dirichlet boundary conditions for the equation
(\ref{HCMA2}) on $M$. Generalized geodesics will correspond to generalized solutions
of the equation (\ref{HCMA2}) in the sense of pluripotential theory, which are invariant under the rotation $w\to e^{i\theta}w$.

\medskip
For the purpose of stability, we are particularly interested in geodesic rays, which correspond to $T=\infty$, and the manifold $M$ is given by $M=X\times D^\times$, with $D^\times=\{0<|w|<1\}$ being the pointed disk. We compactify $M$ into $\hat M=X\times D$, by adjoining the central fiber $X_0=X\times \{0\}$.
Then $\o+ {i\over 2}\ddb |w|^2$ is a \K metric on $M$. Let $p_1,\cdots,p_N\in M$ be any $N$ distinct points
in the central fiber, i.e.,
$\pi(p_\alpha)=0\in\C$, where $\pi:M\to D$
is the projection on the second factor. For each $\alpha$, let
$U_\alpha$ be a neighborhood of $p_\alpha$ in $X$, and
let $f_{1\alpha}(z,w),\cdots, f_{n+1,\alpha}(z,w)$ be any $n+1$
holomorphic functions on $U_\alpha\times D$, with 
the property that their only common zero is at $(p_\alpha,0)$
and with $\sum_{j=1}^n|f_{j\alpha}(z,w)|^2$ invariant under the rotation
$w\to e^{i\theta}w$. Theorem \ref{1} gives then a function
$G(z,w;p_1,\cdots,p_N)$ with the prescribed singularities
(\ref{asymptotics}) and satisfying the equation
$(\o + {i\over 2}\ddb(|w|^2+G))^{n+1}=0$
on $X\times D^\times$. We may choose the blow ups
in the proof of Theorem \ref{1} to be equivariant under the above rotation.
The function $G(z,w;p_1,\cdots,p_N)$ must also be invariant under
rotation, so $|w|^2+G$ defines a generalized
geodesic ray in the space of \K metrics.
Because of the singularities at $(p_\alpha,0)$, the rays are not trivial
(i.e., $\vp$ is not constant along the ray).
They are different from those previously obtained in the literature.

\bigskip
\centerline{}
\bigskip

D.H. Phong

Department of Mathematics, Columbia University, New York, NY 10027

\bigskip

Jacob Sturm

Department of Mathematics, Rutgers University, Newark, NJ 07102

\end{document}